\documentclass{amsart}

\usepackage{amsmath,amsthm,amssymb}
\usepackage{enumerate}
\usepackage{graphicx}
\usepackage{float}

\theoremstyle{plain}
\newtheorem{thm}{Theorem}[section]

\theoremstyle{definition}

\theoremstyle{remark}

\usepackage{tikz}
\usetikzlibrary{arrows,shapes,trees}
\usepackage{graphicx}
\usepackage{hyperref}
\usepackage{amssymb}
\usepackage{amsfonts}
\usepackage{amsthm}
\usepackage{amsmath}
\usepackage{epstopdf}
\usepackage{mathrsfs}

\newcommand{\Z}{\mathbb Z}

\newcommand{\R}{\mathbb R}

\renewcommand{\k}{\mathbf{k}}

\newcommand{\e}{\varepsilon}

\renewcommand{\l}{\ell}
\newcommand{\la}{\lambda}

\title{Discontinuity of Lyapunov exponents near fiber bunched cocycles}
\author{Clark Butler}
\thanks{This material is based upon work supported by the National Science Foundation Graduate Research Fellowship under Grant \# DGE-1144082. }
\begin{document}

\begin{abstract}
We give examples of locally constant $SL(2,\R)$-cocycles over a Bernoulli shift which are discontinuity points for Lyapunov exponents in the H\"older topology and are arbitrarily close to satisfying the fiber bunching inequality. Backes, Brown, and the author have shown that the Lyapunov exponents vary continuously when restricted to the space of fiber bunched H\"older continuous cocycles \cite{BBB}. Our examples give evidence that this theorem is optimal within certain families of H\"older cocycles. 
\end{abstract} 
\maketitle
\section{Introduction}
There is a fundamental dichotomy between the generic behavior of continuous linear cocycles over hyperbolic systems and the generic behavior of cocycles of higher regularity. For our discussion we specialize to the case of cocycles taking values in $SL(2,\R)$. In the continuous setting a theorem of Bochi shows that the generic cocycle over a fixed ergodic dynamical system on a compact space is either uniformly hyperbolic or has all Lyapunov exponents equal to zero \cite{Boc1}; this is proven by showing that these cocycles are the only continuity points for Lyapunov exponents. In sharp contrast, Backes, Brown and the author showed that among H\"older continuous $SL(2,\R)$-valued cocycles which satisfy the \emph{fiber bunching} condition, the Lyapunov exponents vary continuously \cite{BBB}. Fiber bunching is an open condition among H\"older continuous cocycles which is defined below. As a consequence, there are open sets of $SL(2,\R)$-valued H\"older continuous cocycles over hyperbolic systems which are nonuniformly hyperbolic (in the sense that the Lyapunov exponents are nonzero) but are not uniformly hyperbolic.  

A natural question then arises: what is the generic behavior of H\"older continuous cocycles which are not fiber bunched? The strongest theorem currently known in this direction is due to Viana \cite{V08}, who showed that cocycles with nonzero Lyapunov exponents are generic in both the topological and measure-theoretic sense in the space of H\"older  $SL(2,\R)$ cocycles. However Bocker-Neto and Viana \cite{BV} have constructed examples of non-fiber bunched cocycles with nonzero Lyapunov exponents which can be arbitrarily well approximated in the H\"older topology by cocycles with all zero Lyapunov exponents, much like in Bochi's theorem. These examples are very far from being fiber bunched. Viana asked \cite[Ch. 9]{V} whether this construction could be improved to give examples of cocycles which are arbitrarily close to being fiber bunched but are still discontinuity points for Lyapunov exponents. The goal of this paper is to give an affirmative answer to this question. Our construction also raises some related questions regarding the prevalence of discontinuity of Lyapunov exponents away from the fiber bunched cocycles. 

Let $\Sigma = \{0,1\}^{\mathbb{Z}}$ be the space of bi-infinite sequences on two symbols and $f: \Sigma \rightarrow \Sigma$, $f((x_{n})_{n\in \Z}) = (x_{n+1})_{n \in \Z}$ the left shift map.  For $x = (x_{n})_{n \in \mathbb{Z}} \in \Sigma$ we define a distance $d$ by 
\[
d(x,y)=2^{-N(x,y)},\; \textrm{where} \; N(x,y)=\max \lbrace N\geq 0; x_n=y_n \; \textrm{for all} \; |n|<N \rbrace,
\] 
and let $C^{\alpha}(\Sigma,SL(2,\mathbb{R}))$ be the space of $\alpha$-H\"older linear cocycles over $f$ with respect to this metric, equipped with the $\alpha$-H\"older norm defined by (for $A \in C^{\alpha}(\Sigma,SL(2,\mathbb{R}))$)
\[
\|A\|_{\alpha} = \sup_{x\in \Sigma}\|A(x)\| + \sup_{x \neq y \in \Sigma}\frac{\|A(x)-A(y)\|}{d(x,y)^{\alpha}},
\]
where $\| \cdot \|$ is taken to be the standard Euclidean norm on $\R^{2}$ along with the associated operator norm on $SL(2,\R)$. We define $A^{n}(x) = A(f^{n-1}(x)) \cdots A(x)$ and for an $f$-invariant probability measure $\mu$ on $\Sigma$ define the upper and lower Lyapunov exponents of $A$ with respect to $\mu$ to be
\[
\la_{+}(A,\mu) = \inf_{n \geq 1} \frac{1}{n}\int_{\Sigma} \log \|A^{n}\| \, d\mu,
\]
\[
\la_{-}(A,\mu) = \inf_{n \geq 1} \frac{1}{n}\int_{\Sigma} \log \|(A^{n})^{-1}\|^{-1} \, d\mu .
\]
Since $A$ takes values in $SL(2,\R)$, we have the relationship $\la_{-}(A,\mu) = -\la_{+}(A,\mu)$. For $\alpha > 0$, a linear cocycle $A \in C^{\alpha}(\Sigma,SL(2,\R)$ is \emph{$\alpha$-fiber bunched} if there is an $n > 0$ such that 
\[
\sup_{x \in \Sigma}\|A^{n}(x)\| \cdot \|(A^{n}(x))^{-1}\| < 2^{\alpha n}.
\]
In \cite{BBB} it is shown that if $\mu$ is ergodic, fully supported, and has continuous local product structure then the map $A \rightarrow \la_{+}(A,\mu)$ is continuous when restricted to the set of $\alpha$-fiber bunched cocycles in $C^{\alpha}(\Sigma,SL(2,\R)$. 

We define a 1-parameter family of linear cocycles $A_{\sigma}: \Sigma \rightarrow SL(2,\mathbb{R})$ by first defining $A'_{\sigma}$ on the two-point set $\{0,1\}$ by 
\[
A'_{\sigma}(0) = \left[ \begin{array}{cc} \sigma^{-1} & 0 \\
0 & \sigma \end{array}\right], \; A'_{\sigma}(1) = \left[ \begin{array}{cc} \sigma & 0 \\
0 & \sigma^{-1} \end{array}\right],
\]
with $\sigma > 1$. Now let $\pi: \Sigma \rightarrow \{0,1\}$ be projection onto the 0-coordinate and define $A_{\sigma} = A'_{\sigma} \circ \pi$. We also introduce a 1-parameter family of measures given by taking $\nu_{p}$ to be the probability measure on $\{0,1\}$ defined by $\nu_{p}(\{0\}) = 1-p$, $\nu_{p}(\{1\}) = p$ and then defining $\mu_{p}$ on $\Sigma$ as the product measure $\mu_{p} = \nu_{p}^{\mathbb{Z}}$. We will restrict our attention to $p \in (1/2,1)$. It is easy to compute directly from the Birkhoff ergodic theorem that 
\[
\la_{+}(A_{\sigma},\mu_{p}) = (2p-1)\log \sigma. 
\] 
The main theorem of this paper is the following, 
\begin{thm}\label{theorem: discontinuity}
Let $p \in (1/2,1)$. If $\sigma^{4p-2} \geq 2^{\alpha}$ then for each open neighborhood $\mathcal{U} \subset C^{\alpha}(\Sigma,SL(2,\mathbb{R}))$ of $A_{\sigma}$ and every $\kappa \in (0,(2p-1)\log \sigma]$  there is a locally constant cocycle $L \in \mathcal{U}$ such that $\la_{+}(L,\mu_{p}) = \kappa$. In particular $A_{\sigma}$ is a discontinuity point for Lyapunov exponents with respect to $\mu_{p}$ in  $C^{\alpha}(\Sigma,SL(2,\mathbb{R}))$. 
\end{thm}
The cocycles $A_{\sigma}$ satisfy
\[
\sup_{x \in \Sigma} \|A_{\sigma}^{n}(x)\| \cdot \|(A_{\sigma}^{n}(x))^{-1}\|  = \sigma^{2n},
\]
for every $n > 0$. Thus $A_{\sigma}$ is $\alpha$-fiber bunched if and only if $\sigma^{2} < 2^{\alpha}$; in this case $A_{\sigma}$ is a continuity point for Lyapunov exponents with respect to $\mu_{p}$ in $C^{\alpha}(\Sigma, SL(2,\R))$. However, if $\sigma^{2} > 2^{\alpha}$ then we can choose $p \in (1/2,1)$ sufficiently close to 1 such that $\sigma^{4p-2} \geq  2^{\alpha}$ and obtain that $A_{\sigma}$ is a discontinuity point for Lyapunov exponents with respect to $\mu_{p}$ in $C^{\alpha}(\Sigma,SL(2,\R))$. Hence Theorem \ref{theorem: discontinuity} gives a family of examples of discontinuity points for Lyapunov exponents which come arbitrarily close to satisfying the fiber bunching inequality. 

The inequality $\sigma^{4p-2} \geq 2^{\alpha}$ comes from the observation that for $\mu_{p}$-a.e. $x \in \Sigma$, 
\begin{align*}
\lim_{n \rightarrow \infty} \frac{1}{n} \log \left(\|A_{\sigma}^{n}(x)\| \cdot \|(A_{\sigma}^{n}(x))^{-1}\|\right) &= \la_{+}(A_{\sigma},\mu_{p}) -  \la_{-}(A_{\sigma},\mu_{p}) \\
&= (4p-2) \log \sigma .
\end{align*}
Thus if we want the fiber bunching inequality to be violated along a $\mu_{p}$-typical orbit of $f$ we must require $\sigma^{4p-2} \geq 2^{\alpha}$. 

The Bocker-Viana construction shows that $A_{\sigma}$ is a discontinuity point for Lyapunov exponents in $C^{\alpha}(\Sigma,SL(2,\R))$ with respect to $\mu_{p}$ for any $p \in (1/2,1)$ provided that $\sigma^{2} > 2^{4\alpha}$ \cite{BV}.  However, in their example they are able to choose the approximating cocycles $L_{k} \rightarrow A_{\sigma}$ to satisfy $\lambda_{+}(L_{k},\mu_{p}) = 0$ for each $k$, whereas our techniques do not enable us to obtain approximating cocycles with zero Lyapunov exponents. It thus remains an interesting question whether it is possible to construct approximating cocycles with vanishing Lyapunov exponents that come arbitrarily close to satisfying the fiber bunching inequality. 

Another pair of interesting questions arises by considering what is necessary in the behavior of the cocycle $A_{\sigma}$ in order to carry out our construction below. 

\textbf{Question 1}: Is $A_{\sigma}$ a continuity point for Lyapunov exponents $C^{\alpha}(\Sigma,SL(2,\R))$ with respect to $\mu_{p}$ if $\sigma^{2} > 2^{\alpha}$ but $\sigma^{4p-2} < 2^{\alpha}$ and the ratio $\frac{\sigma^{2}}{2^{\alpha}}$ is close to 1? When $\sigma^{2} \geq 2^{4\alpha}$ the Bocker-Viana construction shows that $A_{\sigma}$ is a discontinuity point with respect to $\mu_{p}$ for $p \notin \{0,1/2,1\}$. Proper analysis of our construction shows that this can be improved to $\sigma^{2} \geq 2^{2\alpha}$; in this case one can take $\gamma = 1$ and thus $W = Z$ in the selection of parameters in Section \ref{definitions} and the construction is greatly simplified. However the case of $2^{\alpha} \leq \sigma^{2} < 2^{2\alpha}$ and  $\sigma^{4p-2} < 2^{\alpha}$ remains open. 

\textbf{Question 2}: Define a new cocycle $F_{\sigma}: \Sigma \rightarrow SL(2,\R)$ by setting 
\[
F'_{\sigma}(0) = \left[ \begin{array}{cc} 1 & 0 \\
0 & 1 \end{array}\right], \; F'_{\sigma}(1) = \left[ \begin{array}{cc} \sigma & 0 \\
0 & \sigma^{-1} \end{array}\right],
\]
and then letting $F_{\sigma} = F_{\sigma}' \circ \pi$. Suppose that $\sigma^{2} > 2^{\alpha}$. Is $F_{\sigma}$ a continuity point for Lyapunov exponents in $C^{\alpha}(\Sigma,SL(2,\R))$ with respect to $\mu_{p}$ for some $p \notin \{0,1\}$? The motivation for considering the cocycle $F_{\sigma}$ comes from the observation that no iterate of $F_{\sigma}$ expands the second coordinate in $\R^{2}$, whereas a crucial feature of our construction below is that even though the coordinate in $\R^{2}$ is contracted on average (with respect to $\mu_{p}$) by $A_{\sigma}$, there are still points of $\Sigma$ at which the second coordinate is expanded by the factor $\sigma^{2} > 2^{\alpha}$. 

We thank Aaron Brown and Marcelo Viana for useful discussions regarding this construction. We thank Lucas Backes and Amie Wilkinson for reviewing earlier versions of this paper and offering helpful comments. Lastly we thank the anonymous referee for numerous edits to the preliminary draft of this paper, including crucial corrections to the computations at the end of the paper. 

\section{Selection of Parameters, Definitions}\label{definitions}

We first observe that it suffices to prove Theorem \ref{theorem: discontinuity} in the case that $\sigma^{4p-2} > 2^{\alpha}$. For suppose that Theorem \ref{theorem: discontinuity} holds for $\sigma^{4p-2} > 2^{\alpha}$. Now let $\sigma^{4p - 2} = 2^{\alpha}$, let $\mathcal{U}$ be an open neighborhood of $A_{\sigma}$ in $C^{\alpha}(\Sigma,SL(2,\R))$, and let $\kappa \in (0,(2p-1)\log \sigma]$ be given. Choose $\delta > 0$ small enough that $A_{\sigma+\delta} \in \mathcal{U}$. Then $\kappa \in (0,(2p-1)\log (\sigma+\delta)]$, $(\sigma+\delta)^{4p-2} > 2^{\alpha}$, and $\mathcal{U}$ is an open neighborhood of  $A_{\sigma+\delta}$, so we can find a locally constant cocycle $L \in \mathcal{U}$ with $\la_{+}(L,\mu_{p}) = \kappa$. 

 Thus we assume that $\sigma^{4p-2} > 2^{\alpha}$. Set $A:= A_{\sigma}$, $\mu:=\mu_{p}$. Since $p < 1$ the assumptions of the theorem imply that $\sigma^{2} > 2^{\alpha}$. Let $\kappa \in (0,(2p-1)\log \sigma]$ be given and let $\mathcal{U}$ be a given open neighborhood of $A$ in $C^{\alpha}(\Sigma,SL(2,\R))$. We fix the following parameters in the construction, 
\begin{itemize}
\item $\gamma > 0$ is a rational number chosen small enough that $\sigma^{2} > 2^{(\gamma +1)\alpha}$.
\item $\omega$ is then chosen to be the smallest integer such that $\omega >\gamma^{-1}$. 
\item $\beta > 0$ is chosen small enough that $\sigma^{4p - 2 - 4\beta} > 2^{\alpha}$. 
\end{itemize}

We let $N > 0$ be a large integer such that $\gamma N$ is an integer. Throughout the construction we will use $C$ for any multiplicative constant which is independent of $N$ (though it may depend on the parameters $\gamma, \beta$, etc. chosen above). The value of $C$ may change from line to line within a series of inequalities. For quantities $a$ and $b$ which possibly depend on $N$, we write $a \asymp b$ if there is a constant $C$ independent of $N$ such that $C^{-1}b \leq a \leq Cb$. 

 Set $Z \subset \Sigma$ to be the cylinder 
\[
Z = \{x \in \Sigma: x_{0} = 1,\; x_{i} = 0, \; 1 \leq i \leq N\}.
\]
Let $W \supseteq Z$ be the larger cylinder 
\[
W = \{x \in \Sigma: x_{0} = 1,\; x_{i} = 0, \; 1 \leq i \leq  \gamma N\}.
\]
It is clear that the cylinders $f^{i}(Z)$ are disjoint for $0 \leq i \leq N$ and the same is true for $f^{i}(W)$, $0 \leq i \leq \gamma N$.

For $x \in W$ we define 
\[
\tau_{W}(x) = \inf\{n \geq 1: f^{n}(x) \in W\},
\]
to be the first return time of $x$ to $W$ under the map $f$. We define $f^{\tau_{W}}(x):= f^{\tau_{W}(x)}(x)$, $\tau^{(1)}_{W}:=\tau_{W}$, and then inductively define for $n \geq 1$,  
\[
\tau_{W}^{(n)}(x) = \sum_{j=0}^{n-1} \tau_{W}(f^{\tau_{W}^{(j)}}(x)).
\]
Thus $\tau_{W}^{(n)}(x)$ is the minimal number of iterates necessary for $x$ to return to $W$ exactly $n$ times. We then define
\[
t(x)= \inf\{n \geq 1: f^{\tau_{W}^{(n)}}(x) \in Z\},
\]
to be the first return time of $x$ to $Z$ under the first return map $f^{\tau_{W}}$ to $W$. We also define for $x \in W$ and $n \geq 1$, 
\[
S_{\tau_{W}^{(n)}}(x) =  \sum_{i = 0}^{\tau_{W}^{(n)}(x)-1} \pi (f^{i}(x)).
\]
We similarly define $\tau_{Z}(x)$ to be the first return time of $x$ to $Z$ for $x \in Z$ and $S_{\tau_{Z}}(x) = \sum_{i = 0}^{\tau_{Z}(x)-1} \pi (f^{i}(x))$. We note that $\tau_{Z}(x) = \tau_{W}^{(t(x))}(x)$ by definition. Let $\mu_{Z} = \mu(Z)^{-1}\mu|_{Z}$ be the induced invariant measure for the first return map $f^{\tau_{Z}}$ of $f$ on $Z$. 

We can think of $Z = Z_{N}$ as a family of cylinders in $\Sigma$ parametrized by $N$. As noted above, the first $N$ iterates of $Z$ have no overlaps: $f^{i}(Z) \cap Z = \emptyset$ for $1 \leq i \leq N$.  As a consequence, by a result of Abadi-Vergnes \cite{AbVe}, there is a constant $\eta > 0$ independent of $N$ such that for all $a > 1$ and $N > 0$ we have
\[
\mu_{Z}(\{x \in Z: \tau_{Z}(x) > a \mu(Z)^{-1}\}) \leq \eta e^{-a}.
\] 
We then fix a final parameter $\zeta$, 
\begin{itemize}
\item $\zeta > 1$ is chosen large enough that $ \eta (\log \sigma^{2}) \frac{e^{-\zeta}}{(1-e^{-\zeta})^{2}} < \frac{\kappa}{100}$. 
\end{itemize}

We define a \emph{$W$-return block} (of length $m  \geq \gamma N + 1$) to be a finite string $v = v_{0}v_{1}\dots v_{m-1}$ with $v_{i} \in \{0,1\}$, $v_{0} = 1$, $v_{i} = 0$ for $1 \leq i \leq \gamma N$ and such that there is no subsegment of $v_{\gamma N+1}\dots v_{m-1}$ of the form $1 0 \dots 0$ with $\gamma N$-many 0's.  We define $|v|:= m$ to be the length of the $W$-return block $v$. We denote set of all $W$-return blocks by $\mathscr{V}$. For $v \in \mathscr{V}$ we define
\[
\mathscr{C}(v) = \{x \in \Sigma: x_{i} = v_{i}, \, 0 \leq i \leq |v| - 1, \,  x_{|v|} = 1, \, x_{i} =0, \, |v|+1 \leq i \leq |v|+\gamma N\}, 
\]
to be the cylinder associated to $v$ in $\Sigma$. By construction if $x \in \mathscr{C}(v)$ then $x \in W$ and $\tau_{W}(x) = |v|$. 

We let $e_{1} = (1,0)$ and $e_{2} = (0,1)$ be the standard basis for $\mathbb{R}^{2}$. For $\theta \in \R$ we let
\[
R^{\theta}_{1} = \left[ \begin{array}{cc} 1& \theta \\
0 & 1 \end{array}\right], \; R^{\theta}_{2} = \left[ \begin{array}{cc} 1& 0\\
\theta & 1 \end{array}\right],
\]
be shears by $\theta$ which fix $e_{1}$ and $e_{2}$ respectively. 

\section{Construction}\label{construction}
Let $\e > 0$ be given. We first modify $A$ to create a new locally constant cocycle $B: \Sigma \rightarrow SL(2,\R)$ with $\|A-B\|_{\alpha} < C\e$ and such that for $x \in Z$ the first return map of $B$ to $Z$  satisfies 
\[
B^{\tau_{Z}}(x)(e_{1}) \asymp \e C^{N} \sigma^{-2S_{\tau_{Z}}(x) + \tau_{Z}(x)}e_{2}. 
\]

Define $B_{*}: \Sigma \rightarrow SL(2,\mathbb{R})$ by $B_{*}|_{\Sigma \backslash W} = A$ and $B_{*}|_{W} = A \circ R_{2}^{\e 2^{-\gamma \alpha N}}$. An easy computation shows that $\|A-B_{*}\|_{\alpha} < C\varepsilon$. For each $x \in Z$ we then have 
\[
B_{*}^{N}(x)(e_{1}) = \sigma^{2-N}e_{1} + \sigma^{N-2}\varepsilon 2^{-\gamma \alpha N} e_{2} = \sigma^{N-2}\varepsilon 2^{-\gamma \alpha N} \cdot ( \varepsilon^{-1} 2^{\gamma \alpha N}\sigma^{-2(N-2)}e_{1} + e_{2}).
\]

We choose $N$ large enough that 
\[
\varepsilon^{-1} 2^{\gamma \alpha N} \sigma^{-2(N-2)} < \varepsilon 2^{-\alpha N},
\]
which is possible by our choice of $\gamma$, and then we define $B$ by $B|_{\Sigma \backslash f^{N-1}(Z)} = B_{*}$ and $B|_{f^{N-1}(Z)} = R_{1}^{-\theta} \circ B_{*}$, where $\theta > 0$ small is chosen such that $R_{1}^{-\theta}(B_{*}^{N}(e_{1}))$ is a multiple of $e_{2}$. By our choice of $N$ we can choose $\theta$ to satisfy $\theta \asymp \e 2^{-\alpha N}$ and thus it is again easily verified that $\|B_{*} - B\|_{\alpha} < C\varepsilon$.

The desired assertion about the first return map of $B$ follows from observing that the only modification to $A$ which does not fix the line spanned by $e_{2}$ occurs on the cylinder $f^{N-1}(Z)$. Hence for $x \in Z$, 
\[
B^{N}(x)(e_{1}) \asymp \e \sigma^{N} 2^{-\gamma \alpha N}e_{2},
\]
and thus, observing that $B^{\tau_{Z}}(x)(e_{1}) = A^{\tau_{Z}(x)-N}(f^{N}(x))(B^{N}(x)(e_{1}))$ (since the action of $B$ on $e_{2}$ remains unchanged from $A$ outside of $f^{N-1}(Z)$) and that 
\[
S_{\tau_{Z}}(x) = 1+\sum_{i=N+1}^{\tau_{Z}(x)-1}\pi(f^{i}(x)),
\]
we conclude the desired estimate for $B^{\tau_{Z}}(x)(e_{1})$.

The rest of the construction is devoted to modifying $B$ to create a new locally constant cocycle $L$ satisfying $\|B-L\|_{\alpha} < C\e$ and a similar estimate for the first return map of $L$ on $Z$ moving $e_{2}$ to $e_{1}$, 
\[
L^{\tau_{Z}}(x)(e_{2}) \asymp \sigma^{2S_{\tau_{Z}}(x) -\tau_{Z}(x)}e_{1},
\]
while preserving the estimate for the first return map for $B$. This estimate is significantly harder to obtain; we will only be able to move $e_{2}$ to $e_{1}$ with $L^{\tau_{Z}}$ on a large measure subset of $Z$. The difficulty is that if we are forced to choose $\gamma < 1$ in our parameter selection then the $\mu$-typical  point $x \in Z$ will return to $W$ many times before returning to $Z$ under the orbit of $f$. Once we return to $W$ we shear the vector $e_{1}$ to the left of itself and lose control of how close the iterates of $B$ map $e_{2}$ to $e_{1}$. 

To counteract this issue we first isolate and organize the sequences of $W$-return blocks which have desirable behavior. We define for each $\l \geq 1$ a subset $\mathscr{G}_{\l} \subset \mathscr{V}^{\l}$ to be the set of all $\l$-tuples $(v^{1},\dots,v^{\l})$ of $W$-return blocks for which  
\[
\mathscr{C}(v^{j}) \subset \{x \in W: \tau_{W}(x) < N^{j+1}\},
\]
for $1 \leq j \leq \l - 1$ (if $\l \geq 2$) and 
\[
\mathscr{C}(v^{\l}) \subset \{x \in W: N^{\l+1} \leq \tau_{W}(x) \leq \zeta\mu(Z)^{-1}\} \cap \{x \in W: |S_{\tau_{W}}(x) - p \tau_{W}(x)| \leq \beta \tau_{W}(x)\}.
\]
We will identify an $\l$-tuple $(v^{1},\dots,v^{\l})$ of $W$-return blocks with the concatenated string $v^{1}v^{2}\dots v^{\l}$ of these return blocks and we set $\mathscr{C}(v^{1},\dots,v^{\l}) \subset W$ to be the cylinder associated to this concatenated string which ends in a return to $W$. 

To each point $x \in Z$ we associate the sequence of $W$-return blocks 
\[
(v^{1}(x),\dots,v^{t(x)}(x)) \in \mathscr{V}^{t(x)},
\] 
given by looking at the sequence of returns of $x$ to $W$ before it returns to $Z$. We define the \emph{good set} $G \subset Z$ to be the set of $x \in Z$ whose associated sequence of $W$-return blocks satisfies the following properties,
\begin{itemize}
\item $v^{1}(x) \in \mathscr{G}_{1}$,
\item If $v^{i}(x)$ is a return block with $|v^{i}(x)| \geq N^{2}$ then 
\[
\mathscr{C}(v^{i}(x)) \subset \{x \in W: |S_{\tau_{W}}(x) - p \tau_{W}(x)| \leq \beta \tau_{W}(x)\},
\] 
\item There is no sub-tuple $(v^{i+1}(x),v^{i+2}(x),\dots,v^{i+\omega}(x))$ with $\omega$ members such that 
\[
\mathscr{C}(v^{i+j}(x)) \subset \{x \in W: \tau_{W}(x) < N^{j + 1}\}
\]
for each $1 \leq j \leq \omega $,
\item $\mathscr{C}(v^{t(x)}(x)) \subset \{x \in W: \tau_{W}(x) \geq N^{\omega+1}\}$,
\item $\tau_{Z}(x) \leq \zeta \mu(Z)^{-1}$.
\end{itemize}
The final condition bounding the first return time of $x$ to $Z$ will not be used in this section but will be relevant in Section \ref{analysis}. 

We claim that for each $x \in G$ there is a unique partition of the sequence $(v^{1}(x),\dots,v^{t(x)}(x))$ into subsequences of the form $(v^{i}(x),\dots,v^{i+\l-1}(x))$ which lie in $\mathscr{G}_{\l}$ with $1 \leq \l \leq \omega$. For the existence of such a decomposition, observe first that $v^{1}(x) \in \mathscr{G}_{1}$ by the definition of $G$. Now suppose we have such a partition of the first $k$ terms $(v^{1}(x),\dots,v^{k}(x))$ of the sequence of $W$-return blocks associated to $x$. If $|v^{k+1}(x)| \geq N^{2}$ then the definition of $G$ implies that $v^{k+1}(x) \in \mathscr{G}_{1}$ and so we take the next partition element to be $(v^{k+1}(x))$ on its own. 

If $|v^{k+1}(x)|  < N^{2}$ then there must be some $2 \leq \l \leq \omega$ such that $|v^{k+j}(x)| < N^{j+1}$ for $1 \leq j \leq \l-1$ and $|v^{k+\l}(x)| \geq N^{\l + 1}$. It follows from our definition of $G$ that $(v^{k+1}(x),\dots,v^{k+\l}(x)) \in \mathscr{G}_{\l}$ and we take this to be our next partition element. Repeating this inductively, we obtain a partition of $(v^{1}(x),\dots,v^{t(x)}(x))$ with the desired properties, noting that the requirement $|v^{t(x)}(x)| \geq N^{\omega +1}$ forces the last subsequence to be of the form $(v^{t(x)-\l+1}(x),\dots,v^{t(x)}(x)) \in \mathscr{G}_{\l}$ for some $1 \leq \l \leq \omega$. Uniqueness of the partition constructed can be easily verified, but we will not need it so we omit the proof. 

We can now describe how the final perturbation $L$ will be constructed from $B$. $L$ will be constructed so that for each $x \in G$, on the first $W$-return block $v^{1}(x)$ the matrix $L^{\tau_{W}}$ will map $B^{N}(e_{2})$ onto $e_{1}$.  Then on each successive sequence $(v^{i+1}(x), \dots, v^{i+\l}(x)) \in \mathscr{G}_{\l}$ of the partition of the sequence of $W$-return blocks associated to $x$ constructed above, we will design $L$ so that the $\l$th return map $L^{\tau_{W}^{(\l)}}$ restricted to $\mathscr{C}(v^{i+1}(x),\dots,v^{i+\l}(x))$ leaves the line spanned by $e_{1}$-invariant. The result will be that $L^{\tau_{Z}}(x)$ maps $e_{2}$ to the line spanned by $e_{1}$, as desired. 

We give here a brief explanation of the necessity of the third condition defining $G$ as well as the consideration of multiple $W$-return blocks grouped together, as these requirements introduce substantial complications into the definition of $G$ as well as the construction of $L$. If we only considered single $W$-return blocks then it would be necessary to design the perturbation $L$ such that for each $x \in G$ the first return map of $L$ to $W$ on $\mathscr{C}(v^{i}(x))$ for $2 \leq i \leq t(x)$ leaves the line spanned by $e_{1}$ invariant.  We always have the lower bound $|v^{i}(x)| \geq \gamma N + 1$; we also know by the construction of $B$ that for any $y \in W$ the angle between $B^{\gamma N+1}(y)(e_{1})$ and $e_{1}$ is at least $c > 0$, for some constant $c$ independent of $N$.  In order to make an $\alpha$-H\"older small perturbation that ensures that the first return map to $W$ on $\mathscr{C}(v^{i}(x))$  keeps the line through $e_{1}$ invariant while not affecting the behavior of $B$ outside of $\mathscr{C}(v^{i}(x))$, we must at least wait for an iterate $M > \gamma N+1$ such that the angle between $B^{M}(y)(e_{1})$ and $e_{1}$ is exponentially small in $M$ before performing the shear that returns $B^{M}(y)(e_{1})$ to the line spanned by $e_{1}$. 

Unfortunately, a probabilistic computation shows that if $\gamma < \frac{1}{2}$ then the probability that $|v^{i}(x)| = \gamma N + 2$ (in other words, the probability that some $W$-return block is the minimum possible size) for some $1 \leq i \leq t(x)$ approaches 1 as $N \rightarrow \infty$. This means that when $\gamma$ is small and we consider a points $x \in Z$ returning to $Z$ under the shift map $f$, there will almost surely be at least one $W$-return block which is too short to make the perturbation outlined in the previous paragraph small in the $\alpha$-H\"older norm. Thus it is necessary to consider groups of $W$-return blocks in order to ensure that $\mu_{Z}(G)$ is large. 

For $x \in Z$, $B^{N}(e_{2})$ makes an angle $\asymp \theta$ with $e_{2}$ and lies to the left of $e_{2}$. Let $v \in \mathscr{G}_{1}$ be a good $W$-return block such that $\mathscr{C}(v) \subset Z$. This corresponds to the condition $v = v^{1}(x)$ for some $x \in G$. For $x \in \mathscr{C}(v)$, we have 
\begin{align*}
B^{\tau_{W}}(x)(e_{2}) &= A^{\tau_{W}(x)-N}(B^{N}(x)(e_{2})) \\
&\asymp \sigma^{-N}A^{\tau_{W}(x)-N}(-\theta e_{1} + e_{2}) \\
&\asymp \sigma^{-N}(-\theta \sigma^{2S_{\tau_{W}}(x) - \tau_{W}(x) + N}e_{1} + \sigma^{-2S_{\tau_{W}}(x) + \tau_{W}(x) - N}e_{2}) \\
&= \sigma^{-N}\theta \sigma^{2S_{\tau_{W}}(x) - \tau_{W}(x) + N} (-e_{1} + \theta^{-1}\sigma^{-4S_{\tau_{W}}(x) + 2\tau_{W}(x) - 2N}e_{2}) .
\end{align*}
By the definition of $\mathscr{G}_{1}$ we obtain that for $x \in \mathscr{C}(v)$, 
\[
\sigma^{-4S_{\tau_{W}}(x) +2\tau_{W}(x)} \leq \sigma^{(-4p + 2 + 4\beta)\tau_{W}(x)} .
\]
We claim that for $N$ large enough we have
\[
\theta^{-1}\sigma^{(-4p + 2 + 4\beta)\tau_{W}(x)-2N} < \varepsilon 2^{-\alpha (\tau_{W}(x)+\gamma N)}. 
\]
Put $\xi : = \sigma^{-4p + 2 + 4\beta} 2^{\alpha}$. By the choice of $\beta$ we know that $\xi < 1$. The above inequality is equivalent to requiring 
\[
\xi^{\tau_{W}(x)} < \e \theta \sigma^{2N}2^{-\alpha \gamma N}.
\]
Recalling that $\theta  \asymp \e 2^{-\alpha N}$, this means we must have 
\[
\xi^{\tau_{W}(x)} < C\e^{2} \sigma^{2N}2^{-(1+\gamma)\alpha N},
\]
for some constant $C$ independent of $N$. But $x \in \mathscr{C}(v)$ implies that $\tau_{W}(x) \geq N^{2}$. Hence it suffices to have 
\[
\xi^{N^{2}} < C\e^{2} \sigma^{2N}2^{-(1+\gamma)\alpha N},
\]
which clearly holds for $N$ large enough since $\xi < 1$. 

The cylinder $\mathscr{C}(v)$ has diameter at least $2^{-(|v|+\gamma N)}$, where we recall that $|v| = \tau_{W}(x)$ for any $x \in \mathscr{C}(v)$. There is a small positive number $\delta(v) > 0$ such that the shear $R_{2}^{\delta(v)}(B^{\tau_{W}}(x)(e_{2}))$ is parallel to $e_{1}$. The calculation above shows that we can choose $\delta(v)$ to satisfy $\delta(v) < \e 2^{-(|v|+ \gamma N)}$ for $N$ large enough. We thus define $B_{0} |_{f^{|v|-1}(\mathscr{C}(v))} = R_{2}^{\delta(v)} \circ B$ and on the complement of this collection of cylinders we define $B_{0}$ to coincide with $B$. We then have the inequality $\|B_{0}-B\|_{\alpha} < C \e$.

We now consider those $W$-return blocks $v$ for which $\mathscr{C}(v) \subset W \backslash Z$. This corresponds to considering those $v^{i}(x)$ for which $i \geq 2$ for $x \in G$. We proceed inductively for $1 \leq \l \leq \omega$ and assume that $B_{\l-1}$ has been constructed with $\|B_{\l-1}-B\|_{\alpha} < \l C\e$ and such that on any $W$-return block with $|v| \geq N^{1+\l}$ and $\mathscr{C}(v) \subset W \backslash Z$, the $W$-first return map $B_{\l-1}^{\tau_{W}}$ of $B_{\l-1}$ coincides with $B$. Note this is satisfied for $\l = 1$ because the only modifications on $B_{0}$ are made to the $W$-first return maps on $W$-return blocks $v$ with $\mathscr{C}(v) \subset Z$. At the $\l$th stage of the induction, the cocycle $B_{\l}$ will be obtained from $B_{\l - 1}$ by composition with a matrix of the form $R_{2}^{\delta}$ on some distinguished collection of cylinders in $\Sigma$, with the parameter $\delta$ depending on the cylinder. 

At the $\l$th step we consider those $W$-return blocks for which $N^{1+\l} \leq |v| < N^{2+\l}$ and $v \in \mathscr{G}_{1}$ with the exception that if $\l = \omega$ then we don't impose the upper bound and just consider $v \in \mathscr{G}_{1}$ with $|v| \geq N^{1+\omega}$. For any sequence of $W$-return blocks $(v^{1},\dots,v^{k}) \in \mathscr{G}_{k}$ with $\mathscr{C}(v^{i}) \subset W \backslash Z$ for each $i$, we define $\mathscr{R}(v^{1},\dots, v^{k})$ to be the union of all sets of the form $f^{|w^{1}|+\dots + |w^{j}|}(\mathscr{C}(w^{1},\dots,w^{j},v^{1},\dots,v^{k}))$ with $(w^{1},\dots,w^{j},v^{1},\dots,v^{k}) \in \mathscr{G}_{j+k}$ and $\mathscr{C}(w^{i}) \subset W \backslash Z$ for each $i$. Then we define
\[
\mathscr{C}^{*}(v^{1},\dots,v^{k}) = \mathscr{C}(v^{1},\dots,v^{k}) \backslash \mathscr{R}(v^{1},\dots, v^{k}).
\] 

Let $(u^{1},\dots, u^{k-1},v) \in \mathscr{G}_{k}$ be any good sequence of $W$-return blocks preceding $v$ with $k \leq \l$ (we allow the possibility $k = 1$, for which we only consider $v$ on its own). Then for $x \in \mathscr{C}(u^{1},\dots,u^{k-1},v)$,
\[
B_{\l-1}^{\tau_{W}^{(k)}}(x)(e_{1}) = B^{\tau_{W}}(f^{\tau_{W}^{(k-1)}}(x))(B_{\l-1}^{\tau_{W}^{(k-1)}}(x)(e_{1})).
\]

Further we know that $B_{\l-1}^{\tau_{W}^{(k-1)}}(x)(e_{1})$ has the form 
\[
B_{\l-1}^{\tau_{W}^{(k-1)}}(x)(e_{1}) \asymp \sigma^{2 S_{\tau_{W}^{(k-1)}}(x) - \tau_{W}^{(k-1)}(x)}(e_{1} + r(x) e_{2}),
\]
with $|r(x)| \leq C\sigma^{\sum_{j=1}^{\k-1}N^{j+1}}$ (if $k = 1$ then we do not have to worry about this). This bound is derived from the observation that $(u^{1},\dots,u^{k-1},v) \in \mathscr{G}_{k}$ implies that $|u^{j}| \leq N^{j+1}$ for each $1 \leq j \leq k-1$. Then 
\[
B_{\l-1}^{\tau_{W}^{(k)}}(x)(e_{1}) \asymp \sigma^{2 S_{\tau_{W}^{(k)}}(x) - \tau_{W}^{(k)}(x)}(e_{1} + r(x) \sigma^{-4S_{\tau_{W}}(f^{\tau_{W}^{(k-1)}}(x)) +2\tau_{W}(f^{\tau_{W}^{(k-1)}}(x))}\varepsilon 2^{-\gamma \alpha N}e_{2}).
\]

As in the case of the $B_{0}$ modification above we claim there is a small positive number $\delta(u^{1},\dots,u^{k-1},v)$ such that $R^{-\delta(u^{1},\dots,u^{k-1},v)}_{2}(B^{\tau_{W}^{(k)}}(x)(e_{1}))$ is parallel to $e_{1}$. In fact we can see above that we should take 
\[
\delta(u^{1},\dots, u^{k-1},v) \asymp  r(x) \sigma^{-4S_{\tau_{W}}(f^{\tau_{W}^{(k-1)}}(x)) +2\tau_{W}(f^{\tau_{W}^{(k-1)}}(x))}\varepsilon 2^{-\gamma \alpha N}.
\]
We will show that for $N$ large enough,
\[
r(x) \sigma^{-4S_{\tau_{W}}(f^{\tau_{W}^{(k-1)}}(x)) +2\tau_{W}(f^{\tau_{W}^{(k-1)}}(x))}\varepsilon 2^{-\gamma \alpha N} < \e 2^{-\alpha(\sum_{j = 1}^{\l-1} N^{j+1}+\tau_{W}(f^{\tau_{W}^{(k-1)}}(x))+\gamma N)}.
\]
Rearranging and using our bound on $r(x)$ we see that it suffices to have
\[
\sigma^{-4S_{\tau_{W}}(f^{\tau_{W}^{(k-1)}}(x)) +2\tau_{W}(f^{\tau_{W}^{(k-1)}}(x))}2^{\alpha \tau_{W}(f^{\tau_{W}^{(k-1)}}(x))} \leq C 2^{-\alpha \sum_{j = 1}^{\l-1} N^{j+1}} \sigma^{-\sum_{j=1}^{\k-1}N^{j+1}}.
\]
Noting that $\tau_{W}(f^{\tau_{W}^{(k-1)}}(x))\geq N^{1+\l}$ by our assumption on $x$ and recalling that $\xi = \sigma^{-4p + 2 + 4\beta} 2^{\alpha} < 1$, it suffices to show that for $N$ large enough the inequality
\[
\xi^{N^{\l+1}} < C (\sigma^{-1} 2^{-\alpha})^{\sum_{j=1}^{\l-1}N^{j+1}},
\]
holds. For each fixed $\l$ this inequality is true for $N$ large enough: after taking logarithms on each side and rearranging this is equivalent to the inequality
\[
N^{\l + 1} > \frac{1}{\log(\xi^{-1})}(\log C + \log(\sigma 2^{\alpha})\sum_{j=1}^{\l-1}N^{j+1}),
\]
and the right side is a polynomial in $N$ of degree $\l$ which is thus eventually dominated by $N^{\l + 1}$. We choose $N$ large enough that the desired inequality holds for $1 \leq \l \leq \omega$. We define a new cocycle $B_{\l}$ by 
\[
B_{\l}|_{f^{|u^{1}|+\dots + |u^{k-1}|+ |v|-1}(\mathscr{C}^{*}(u^{1},\dots,u^{k-1},v))} = R_{2}^{\delta(u^{1},\dots,u^{k-1},v)} \circ B_{\l-1},
\]
 and let $B_{\l}$ coincide with $B_{\l-1}$ outside any of the cylinders of the above form. We then have $B_{\l}^{\tau_{W}^{(k)}}(x)(e_{1}) \asymp \sigma^{2S_{\tau_{W}^{(k)}}(x) -\tau_{W}^{(k)}(x)}e_{1}$ for $x$ in one of these modification sets $\mathscr{C}^{*}(u^{1},\dots,u^{k-1},v)$. Each set $\mathscr{C}^{*}(u^{1},\dots,u^{k-1},v)$ is a union of cylinders of diameter at least $2^{-\sum_{j=1}^{\l-1}N^{j+1} -|v| - \gamma N}$, so that the inequality $\xi^{N^{\l+1}} < C (\sigma 2^{-\alpha})^{-\sum_{j=1}^{\l-1}N^{j+1}}$ together with the induction hypothesis $\|B_{\l-1}-B\|_{\alpha} < \l C \e$ implies that$\|B_{\l}-B\|_{\alpha} \lesssim (\l + 1)C\e$. Note that we perform these modifications on $\mathscr{C}^{*}(u^{1},\dots,u^{k-1},v)$ instead of  $\mathscr{C}(u^{1},\dots,u^{k-1},v)$ so that all of the cylinders we modify $B_{\l-1}$ on are actually disjoint. 

We set $L := B_{\omega}$. Since all of the modifications from $B$ to $L$ consisted of shears which fixed $e_{2}$, we have the exact same estimate for $L^{\tau_{Z}}(x)(e_{1})$ with $x \in Z$, 
\[
L^{\tau_{Z}}(x)(e_{1}) \asymp \e C^{N}\sigma^{-2S_{\tau_{Z}}(x) + \tau_{Z}(x)}e_{2}.
\]
For $x \in G$ we claim that we also have 
\[
L^{\tau_{Z}}(x)(e_{2})  \asymp \sigma^{2S_{\tau_{Z}}(x) - \tau_{Z}(x)}e_{1}.
\]

To see this, consider the sequence of $W$-return blocks $(v^{1}(x),\dots,v^{t(x)}(x))$ associated to the return of $x$ to $Z$. We have $L^{\tau_{W}}(x) = B_{0}^{\tau_{W}}(x)$ and consequently 
\[
L^{\tau_{W}}(x)(e_{2}) \asymp \sigma^{2S_{\tau_{W}}(x) - \tau_{W}(x)}e_{1}.
\]
We then have a unique partition of $(v^{1}(x),\dots,v^{t(x)}(x))$ into subsequences of the form $(v^{i}(x),\dots,v^{i+\l-1}(x))$ which lie in $\mathscr{G}_{\l}$ with $1 \leq \l \leq \omega$. On $\mathscr{C}(v^{i}(x),\dots,v^{i+\l-1}(x))$ we see that $L^{\tau_{W}^{(\l)}}$ coincides with $B^{\tau_{W}^{(\l)}}_{\l}$, and using this fact we easily show inductively that 
\[
L^{\tau_{Z}}(x)(e_{2}) \asymp \sigma^{-2S_{\tau_{Z}}(x) + \tau_{Z}(x)}e_{1},
\]
as desired. 

We make some final remarks before moving on to the analysis of the Lyapunov exponents of $L$ in Section \ref{analysis}. By construction $L$ is constant on all cylinders of diameter at most $2^{-\zeta\mu(Z)^{-1}-\gamma N}$, hence $L$ is locally constant and in particular is $\alpha$-H\"older continuous. Since $L$ is constructed from $A$ by composing $A$ with different shearing matrices close to the identity on this finite collection of cylinders of diameter at least $2^{-\zeta\mu(Z)^{-1}-N}$, it is easy to see that there is a continuous family of cocycles $L_{t}$, $t \in [0,1]$, such that $L_{0} = A$, $L_{1} = L$, and for $t \in [0,1]$ we have that $L_{t}$ is constant on any cylinder of diameter at most $2^{-\zeta\mu(Z)^{-1}-\gamma N}$, i.e., the cocycles $L_{t}$ are all locally constant on the same collection of cylinders. We also have $\|A-L_{t}\|_{\alpha} < \omega C\e$ for each $t \in [0,1]$. We can thus choose $\e$ small enough that $L_{t}$ lies in the given neighborhood $\mathcal{U}$ for each $t \in [0,1]$. We also choose $\e$ small enough that $\|L_{t}\| \leq \sigma^{2}$ for each $t \in [0,1]$. We will use this bound for $\|L\|$ in Section \ref{analysis}. 

\section{Analysis of the Lyapunov exponents of L}\label{analysis}
We show in this section that for any given $0 < \kappa \leq (2p-1) \log \sigma$ we can find $N$ large enough that we have $\lambda_{+}(L,\mu) < \kappa$. Using the observation at the end of Section \ref{construction}, there is a continuous family of cocycles $\{L_{t}\}_{t \in [0,1]}$ which are locally constant on the same family of cylinders with $L_{0} = A$ and $L_{1} = L$. By the main theorem of \cite{BBB} the map $t \rightarrow \la_{+}(L_{t},\mu)$ is continuous and hence surjects onto $[\la_{+}(L,\mu),\log(\sigma^{2p-1})]$. Since $\la_{+}(L,\mu) < \kappa$, this implies that there is some $t \in [0,1]$ such that $\la_{+}(L_{t},\mu) = \kappa$ which then completes the proof of Theorem \ref{theorem: discontinuity}. 

It remains to show that for $N$ large enough we have $\lambda_{+}(L,\mu) < \kappa$. $\lambda_{+}(L^{\tau_{Z}},\mu_{Z})$ is related to $\lambda_{+}(L,\mu)$ by the formula
\[
\lambda_{+}(L^{\tau_{Z}},\mu_{Z}) = \mu(Z)^{-1}\lambda_{+}(L,\mu),
\] 
and we can bound $\lambda_{+}(L^{\tau_{Z}},\mu_{Z})$ using the $k$th return map to $Z$ by 
\[
\lambda_{+}(L^{\tau_{Z}},\mu_{Z}) \leq \frac{1}{2}\int_{Z}\log \left\|L^{\tau_{Z}^{(2)}}\right\|\,d\mu_{Z} 
\]
(see \cite{V}). We will show for $N$ large enough that 
\[
\frac{1}{2}\int_{Z}\log\left\|L^{\tau_{Z}^{(2)}}\right\| \, d\mu_{Z} < \kappa \mu(Z)^{-1} .
\]
This immediately implies that $\la_{+}(L,\mu) < \kappa$. Proving this assertion will occupy the rest of this section. 

We first analyze the contribution of those $x$ for which $\tau_{Z}(x) > \zeta \mu(Z)^{-1}$ to the integral above. For $a \geq 1$ we recall from Section \ref{definitions} that 
\[
Q_{a} = \{x \in Z: \tau_{Z}(x) > a \cdot \mu(Z)^{-1}\},
\]
and that we have $\mu_{Z}(Q_{a}) \leq \eta e^{-a}$ for $a > 1$. We then have 
\begin{align*}
\int_{Q_{\zeta}} \log \left\|L^{\tau_{Z}}\right\| \,d\mu_{Z} &= \sum_{n=1}^{\infty} \int_{Q_{n \zeta} \backslash Q_{(n+1)\zeta}} \log \left\|L^{\tau_{Z}}\right\| \,d\mu_{Z} \\
&\leq \sum_{n =1}^{\infty} (n+1)\zeta \mu(Z)^{-1} (\log \sigma^{2})\mu_{Z}(Q_{n\zeta}) \\
&\leq \eta  (\log \sigma^{2}) \mu(Z)^{-1} \sum_{n=1}^{\infty}(n+1)\zeta e^{-n\zeta} \\
&\leq \eta (\log \sigma^{2}) \frac{e^{-\zeta}}{(1-e^{-\zeta})^{2}} \mu(Z)^{-1} .
\end{align*}
By our choice of the parameter $\zeta$, this implies that 
\[
\int_{Q_{\zeta}} \log \left\|L^{\tau_{Z}}\right\| \,d\mu_{Z} < \frac{\kappa}{100}\mu(Z)^{-1}.
\]
We also note that our choice of $\zeta$ implies $\mu_{Z}(Q_{\zeta}) < \frac{\kappa}{100\log \sigma^{2}}$. 
 
We define $K = Z\backslash (G \cup Q_{\zeta})$ and next estimate $\mu_{Z}(K)$. There are four ways for a point $x \in Z$ to fail to belong to $G$ once we exclude the long returns $\tau_{Z}(x) > \zeta \mu(Z)^{-1}$, 
\begin{enumerate}
\item $|v^{1}(x)| < N^{2}$,

\item There is some $1 \leq i \leq t(x)$ such that $|v^{i}(x)| \geq N^{2}$ but 
\[
\left|S_{\tau_{W}}(f^{\tau_{W}^{(i-1)}}(x)) - p \tau_{W}(f^{\tau_{W}^{(i-1)}}(x))\right| > \beta \tau_{W}(f^{\tau_{W}^{(i-1)}}(x)),
\] 

\item There is some $1 \leq i \leq t(x) - \omega $ such that for each $1 \leq j \leq \omega$ we have $|v^{i+j-1}(x)| < N^{j+1}$,

\item $|v^{t(x)}(x)| < N^{\omega+1}$.
\end{enumerate}
We deal with condition (2) first. The measure of the set of $x \in Z$ which satisfy (2) is bounded above by 
\[
\sum_{i=1}^{\zeta \mu(Z)^{-1}}\mu_{Z}\left(\left\{x \in Z: \tau_{W}(f^{\tau_{W}^{(i-1)}}(x)) \geq N^{2}, \,\left|S_{\tau_{W}}(f^{\tau_{W}^{(i-1)}}(x)) - p \tau_{W}(f^{\tau_{W}^{(i-1)}}(x))\right| > \beta \tau_{W}(f^{\tau_{W}^{(i-1)}}(x))\right\}\right).
\]
By summing over the discrete set of possible values of $\tau_{W}(f^{\tau_{W}^{(i-1)}}(x))$, this is bounded above by 
\[
\sum_{i=1}^{\zeta \mu(Z)^{-1}}\sum_{n=N^{2}}^{\infty}\mu_{Z}\left(\left\{x \in Z:  \left|S_{n}(f^{\tau_{W}^{(i-1)}}(x)) - p n \right| > \beta n\right\}\right).
\]
We note that if we define $\mu_{W} = \mu(W)^{-1}\mu|_{W}$, then $\mu_{W}$ is likewise invariant under the first return map $f^{\tau_{W}}$ to $W$ and the above sum is bounded above by 
\begin{align*}
&\frac{\mu(W)}{\mu(Z)}\sum_{i=1}^{\zeta \mu(Z)^{-1}}\sum_{n=N^{2}}^{\infty}\mu_{W}\left(\left\{x \in W:  \left|S_{n}(f^{\tau_{W}^{(i-1)}}(x)) - p n \right| > \beta n\right\}\right) 
\\ &= \frac{\zeta\mu(W)}{\mu(Z)^{2}}\sum_{n=N^{2}}^{\infty}\mu_{W}\left(\left\{x \in W:  \left|S_{n}(x) - p n \right| > \beta n\right\}\right).
\end{align*}
By the well-known Chernoff inequality for sums of bounded, independent, identically distributed random variables (for a quick proof we refer to \cite[Ch. 1]{T12}), 
\[
\mu_{W}\left(\left\{x \in W:  \left|S_{n}(x) - p n \right| > \beta n\right\}\right) \leq e^{-\frac{n\beta^{2}}{2}}.
\]
Hence we at last conclude that 
\begin{align*}
\frac{\zeta\mu(W)}{\mu(Z)^{2}}\sum_{n=N^{2}}^{\infty}\mu_{W}\left(\left\{x \in W:  \left|S_{n}(x) - p n \right| > \beta n\right\}\right) &\leq \frac{\zeta\mu(W)}{\mu(Z)^{2}}\sum_{n=N^{2}}^{\infty}e^{-\frac{n\beta^{2}}{2}} \\
&= \frac{\zeta\mu(W)}{\mu(Z)^{2}}\frac{e^{-\frac{N^{2}\beta^{2}}{2}}}{1-e^{-\frac{\beta^{2}}{2}}}.
\end{align*}

For the other three conditions, we will need a crude estimate on the probability of short return times to $W$ and short hitting times from $Z$ to $W$.  We first observe that for $x = (x_{n})_{n \in \Z} \in Z$, the coordinates $x_{n}$ with $n \geq N + 1$ or $n < 0$ are independent and identically distributed with respect to $\mu_{Z}$, with $\mu_{Z}(\{x_{n} = 1\}) = p$ and $\mu_{Z}(\{x_{n} = 0\}) = 1-p$ for $n \geq N+1$ or $n < 0$. Our second observation is that if $x \in Z$ and there is an $m > 0$ such that $x_{m} = 1$ and $x_{m+i} = 0$ for $1 \leq i \leq \gamma N$ (this corresponds to $f^{m}(x) \in W$) then in fact $m \geq N$. Hence for each $m > N$ we have
\begin{align*}
\mu_{Z}(\{x \in Z: \tau_{W}(x) = m\}) &\leq \mu_{Z}(\{x \in Z: f^{m}(x) \in W\}) \\
&= \mu_{Z}(\{x_{m} = 1\}) \prod_{i=1}^{\gamma N}\mu_{Z}(\{x_{m+i} = 0\}) \\
&= p(1-p)^{\gamma N} = \mu(W).
\end{align*}
Consequently we have 
\begin{align*}
\mu_{Z}(\{x \in Z: \tau_{W}(x) <  m\}) &= \sum_{j = N+1}^{m-1} \mu_{Z}(\{x \in Z: \tau_{W}(x) = j\}) \\
&\leq m \mu(W).
\end{align*}
For condition (1) we thus obtain
\[
\mu_{Z}(\{x \in Z: |v^{1}(x)| < N^{2}\}) \leq N^{2}\mu(W).
\]
For condition (4) we use the invariance of $\mu_{Z}$ under $f^{\tau_{Z}}$ to obtain 
\begin{align*}
\mu_{Z}(\{x \in Z: |v^{t(x)}(x)| < N^{\omega + 1}\}) &\leq \sum_{m=1}^{N^{\omega+1}}\mu_{Z}(\{x \in Z: f^{-m}(f^{\tau_{Z}}(x)) \in W\}) \\
&= \sum_{m=1}^{N^{\omega+1}}\mu_{Z}(\{x \in Z: f^{-m}(x) \in W\}).
\end{align*}

Using the independence of the coordinates $x_{n}$ of $x \in Z$ under $\mu_{Z}$ for $n < 0$ we then see, similarly to the estimates we used for (1), that 
\begin{align*}
\sum_{m=1}^{N^{\omega+1}}& \mu_{Z}(\{x \in Z: f^{-m}(x) \in W\})  \\
&\leq \sum_{m=1}^{N^{\omega+1}}\mu_{Z}(\{x_{-m} = 1\}) \prod_{i=1}^{\gamma N}\mu_{Z}(\{x_{-m+i} = 0\}) \leq N^{\omega +1}\mu(W).
\end{align*}

For condition (3) we can again directly estimate the probability of a $W$-return block sequence $(v^{i}(x),\dots,v^{i+ \omega -1}(x))$ occurring with $|v^{i+j-1}(x)| < N^{j+1}$ for $1 \leq j \leq \omega$. Using the independence of the coordinates $x_{n}$ of $x$ for $n > N$, we have the bound, for each $i \geq 1$, 
\begin{align*}
\mu_{Z}(\{x &\in Z: \tau_{W}^{(i+j-1)}(x) -  \tau_{W}^{(i+j-2)}(x) < N^{j+1},\,1\leq j \leq \omega \})  \\ = &\prod_{j =1}^{\omega} \mu_{Z}(\{x \in Z: \tau_{W}^{(i+j-1)}(x) - \tau_{W}^{(i+j-2)}(x) < N^{j+1} \}) \\
&\leq \prod_{j =1}^{\omega} N^{j+1}\mu(W) \\
& \leq N^{\frac{(\omega +1)(\omega +2)}{2} }\mu(W)^{\omega}.
\end{align*}

We must sum this bound over all possible values of $i$. Because we require $\tau_{Z}(x) \leq \zeta \mu(Z)^{-1}$, we can sum over $1 \leq i \leq \zeta \mu(Z)^{-1}$ and thus bound the $\mu_{Z}$-measure of the set of $x \in Z$ satisfying condition (3) by $\zeta \mu(Z)^{-1}N^{\frac{(\omega +1)(\omega +2)}{2} }\mu(W)^{\omega}$.

Combining our estimates for conditions (1)-(4), we thus conclude that 
\[
\mu_{Z}(K) \leq N^{2}\mu(W) + \frac{\zeta\mu(W)}{\mu(Z)^{2}}\frac{e^{-\frac{N^{2}\beta^{2}}{2}}}{1-e^{-\frac{\beta^{2}}{2}}} + \zeta \mu(Z)^{-1}N^{\frac{(\omega +1)(\omega +2)}{2} }\mu(W)^{\omega} + N^{\omega +1}\mu(W).
\]
Recall that $\mu(W) = p(1-p)^{\gamma N}$ and $\mu(Z) = p(1-p)^{N}$. Plugging these values in for the measures of these sets, we get 
\begin{align*}
\mu_{Z}(K) &\leq N^{2}p(1-p)^{\gamma N} + \zeta p^{-1}(1-p)^{(\gamma-2)N}\frac{e^{-\frac{N^{2}\beta^{2}}{2}}}{1-e^{-\frac{\beta^{2}}{2}}} \\
&+ \zeta N^{\frac{(\omega +1)(\omega +2)}{2} }p^{\omega-1}(1-p)^{(\omega \gamma - 1)N}+ N^{\omega +1}p(1-p)^{\gamma N}.
\end{align*}
Recall that we chose $\omega$ such that $\omega \gamma > 1$. Hence we conclude that as $N \rightarrow \infty$, the right side converges to 0 as $N \rightarrow \infty$. Choose $N$ large enough that 
\[
\mu_{Z}(K) < \frac{\kappa }{100\zeta \log \sigma^{2} }.
\]
We conclude that 
\[
\mu_{Z}(Z \backslash G) < \frac{\kappa }{50\zeta \log \sigma^{2} },
\]
and we also have 
\begin{align*}
\int_{Z \backslash G} \log \left\|L^{\tau_{Z}}\right\| \,d\mu_{Z} &= \int_{K} \log \left\|L^{\tau_{Z}}\right\| \,d\mu_{Z} + \int_{Q_{\zeta}} \log \left\|L^{\tau_{Z}}\right\| \,d\mu_{Z} \\
&\leq \zeta \mu(Z)^{-1}\mu_{Z}(K) \log \sigma^{2} + \int_{Q_{\zeta}} \log \left\|L^{\tau_{Z}}\right\| \,d\mu_{Z}\\
&< \frac{\kappa}{50}\mu(Z)^{-1}.
\end{align*}

We now use the bound 
\[
\frac{1}{2}\int_{Z \backslash G} \log \left\|L^{\tau_{Z}^{(2)}}\right\| \,d\mu_{Z} \leq \frac{1}{2}\int_{Z \backslash G} \log \left\|L^{\tau_{Z}}\right\| \,d\mu_{Z} + \frac{1}{2}\int_{Z \backslash G} \log \left\|L^{\tau_{Z}} \circ f^{\tau_{Z}}\right\| \,d\mu_{Z}.
\]
The first integral is bounded by $\frac{\kappa}{100}\mu(Z)^{-1}$. For the second integral, we observe that the characteristic function of $Z \backslash G$ and the function $L^{\tau_{Z}} \circ f^{\tau_{Z}}$ are independent with respect to $\mu_{Z}$ and thus
\begin{align*}
\frac{1}{2}\int_{Z \backslash G} \log \left\|L^{\tau_{Z}} \circ f^{\tau_{Z}}\right\| \,d\mu_{Z} &= \frac{1}{2}\mu_{Z}(Z \backslash G) \int_{Z} \log \left\|L^{\tau_{Z}} \circ f^{\tau_{Z}}\right\| \,d\mu_{Z} \\
&\leq \frac{1}{2}\mu_{Z}(Z \backslash G) (\log \sigma^{2}) \int_{Z} \tau_{Z} \circ f^{\tau_{Z}} \,d\mu_{Z} \\
&= \frac{1}{2}\mu_{Z}(Z \backslash G) (\log \sigma^{2})\mu(Z)^{-1} \\
&< \frac{\kappa}{100}\mu(Z)^{-1}.
\end{align*}
Thus we conclude that
\[
\frac{1}{2}\int_{Z \backslash G} \log \left\|L^{\tau_{Z}^{(2)}}\right\| \,d\mu_{Z} < \frac{\kappa}{50}\mu(Z)^{-1}.
\]
Applying the exact same calculations to $Z \backslash f^{-\tau_{Z}}(G)$, we obtain also that 
\[
\frac{1}{2}\int_{Z \backslash f^{-\tau_{Z}}(G)} \log \left\|L^{\tau_{Z}^{(2)}}\right\| \,d\mu_{Z} < \frac{\kappa}{50}\mu(Z)^{-1}.
\]
We conclude that 
\[
\frac{1}{2}\int_{Z \backslash (G \cap f^{-\tau_{Z}}(G))}\log \left\|L^{\tau_{Z}^{(2)}}\right\| \,d\mu_{Z} < \frac{\kappa}{25}\mu(Z)^{-1}.
\]

For $x \in G \cap f^{-\tau_{Z}}(G)$, the conclusion of Section \ref{construction} implies that 
\[
\log \left\|L^{\tau_{Z}^{(2)}}(x)\right\| \leq N\log C + (\log \sigma^{2}) \left|2S_{\tau_{Z}}(x) - \tau_{Z}(x) -  2S_{\tau_{Z}}(f^{\tau_{Z}}(x)) + \tau_{Z}(f^{\tau_{Z}}(x))\right| .
\]
 
Define for $x \in Z$, 
\[
T_{n}(x) = \sum_{j=N+1}^{n-1} (2\pi(f^{j}(x)) - 1) = 2S_{n}(x) - n + N .
\]
We view $T_{n}$  as a sum of the independent, identically distributed (with respect to $\mu_{Z}$) random variables $2(\pi \circ f^{j}) -1$ for $j \geq N+1$ and note that $\tau_{Z}-N$ is a stopping time with respect to the natural filtration of the Borel $\sigma$-algebra of $Z$ induced by this sequence of random variables. For each $j \geq N+1$ we have 
\[
\int_{Z} (2(\pi \circ f^{j}) -1) \,d\mu_{Z} = 2p-1
\]
and
\[
\int_{Z} (2(\pi \circ f^{j})-1 - (2p-1))^{2} \,d\mu_{Z} = 1- (2p-1)^{2} .
\]

The Wald identities for a sequence of independent, identically distributed random variables $X_{1},\dots,X_{n}$ with finite second moment states that if $\tau$ is a stopping time with respect to the sequence of $\sigma$-algebras associated with this sequence of random variables, $S_{n} = \sum_{i=1}^{n} X_{i}$, and $S_{\tau}$ is the stopped sum $S_{\tau} = \sum_{i=1}^{\tau} X_{i}$, then the expectations of these random variables satisfy
\[
\mathbb{E}[S_{\tau}] = \mathbb{E}[X_{1}] \cdot \mathbb{E}[\tau],
\]
\[
\mathbb{E}[(S_{\tau}-\mathbb{E}[S_{\tau}])^{2}] = \mathbb{E}[(X_{1}-\mathbb{E}[X_{1}])^{2}] \cdot \mathbb{E}[\tau]. 
\]

We apply these identities to $T_{n}$ stopped at $\tau_{Z}-N$ to conclude that 
\[
\int_{Z}T_{\tau_{Z}}\, d\mu_{Z} = (2p - 1)(\mu(Z)^{-1}-N),
\]
by the first identity and thus by the second identity, 
\[
\int_{Z}(T_{\tau_{Z}}-(2p - 1)(\mu(Z)^{-1}-N))^{2}\, d\mu_{Z} = (1- (2p-1)^{2}) (\mu(Z)^{-1}-N).
\]

Define 
\begin{align*}
\psi(x) &= 2S_{\tau_{Z}}(x) - \tau_{Z}(x) -  2S_{\tau_{Z}}(f^{\tau_{Z}}(x)) + \tau_{Z}(f^{\tau_{Z}}(x)) \\
&= T_{\tau_{Z}}(x) - T_{\tau_{Z}}(f^{\tau_{Z}}(x))
\end{align*}
We clearly have $\int_{Z} \psi \, d\mu_{Z} = 0$. Using the independence of the random variables $T_{\tau_{Z}}$ and $T_{\tau_{Z}} \circ f^{\tau_{Z}}$ together with the fact that they are identically distributed, we conclude that 
\begin{align*}
\int_{Z}\psi^{2}\,d\mu_{Z} &= 2\int_{Z}(T_{\tau_{Z}} - (2p - 1)(\mu(Z)^{-1}-N))^{2}d\mu_{Z}  \\
&= 2(1- (2p-1)^{2} ) (\mu(Z)^{-1}-N) \\
&\leq C \mu(Z)^{-1},
\end{align*}
for $N$ large. Thus by Chebyshev's inequality, for any $a > 0$ we have 
\[
\mu_{Z}(\{x \in Z: |\psi(x)| > a\}) \leq \frac{C}{\mu(Z)a^{2}}.
\]
Putting $a = \mu(Z)^{-3/4}$, this gives
\[
\mu_{Z}(\{x \in Z: |\psi(x)| > \mu(Z)^{-3/4}\}) \leq C\mu(Z)^{1/2}.
\]
Hence 
\begin{align*}
\int_{G \cap f^{-\tau_{Z}}(G)} |\psi | \,d\mu_{Z} &\leq C\mu(Z)^{-\frac{3}{4}} + 2\zeta \mu(Z)^{-1}\mu(Z)^{1/2} \\
&\leq C  \mu(Z)^{-3/4},
\end{align*}
where on the set on which $|\psi| > \mu(Z)^{-3/4}$ we used the upper bound $\tau_{Z}^{(2)}(x) \leq 2\zeta \mu(Z)^{-1}$ imposed on the second return time by membership in $G \cap f^{-\tau_{Z}}(G)$. We choose $N$ large enough that $C \mu(Z)^{-3/4} < \frac{\kappa}{25}\mu(Z)^{-1}$. We then finally get the conclusion 
\[
\frac{1}{2}\int_{Z}\log \left\|L^{\tau_{Z}^{(2)}}\right\| \,d\mu_{Z}  <  \frac{2\kappa}{25} \mu(Z)^{-1} + N \log C \leq \kappa \mu(Z)^{-1},
\]
for $N$ large enough, which completes the proof. 
\bibliographystyle{plain}
\bibliography{lyapdiscontinuity}
\end{document}